\pdfoutput=1
\RequirePackage{ifpdf}
\ifpdf 
\documentclass[pdftex]{sigma}
\else
\documentclass{sigma}
\fi

\numberwithin{equation}{section}

\newcommand{\eqdef}{\stackrel{\text{def}}{=}}
\newcommand{\n}{\nonumber\\}

\newcommand{\ignore}[1]{}

\newcommand{\Romannumeral}[1]{\uppercase\expandafter{\romannumeral#1}}

\newcommand{\ai}{\text{I}}
\newcommand{\ait}{\text{II}}
\newcommand{\aitr}{\text{III}}

\begin{document}


\renewcommand{\thefootnote}{$\star$}

\renewcommand{\PaperNumber}{085}

\FirstPageHeading

\ShortArticleName{Global Solutions  with
                         a High Degree of Apparent Singularity}

\ArticleName{Global Solutions of Certain Second-Order\\ Dif\/ferential Equations with
                         a High Degree\\ of Apparent Singularity\footnote{This
paper is a contribution to the Special Issue ``Superintegrability, Exact Solvability, and Special Functions''. The full collection is available at \href{http://www.emis.de/journals/SIGMA/SESSF2012.html}{http://www.emis.de/journals/SIGMA/SESSF2012.html}}}

\Author{Ryu SASAKI~$^\dag$ and Kouichi TAKEMURA~$^\ddag$}

\AuthorNameForHeading{R.~Sasaki and K.~Takemura}

\Address{$^\dag$~Yukawa Institute for Theoretical Physics,
     Kyoto University, Kyoto 606-8502, Japan}
\EmailD{\href{mailto:ryu@yukawa.kyoto-u.ac.jp}{ryu@yukawa.kyoto-u.ac.jp}}

\Address{$^\ddag$~Department of Mathematics, Faculty of Science and Technology,\
      Chuo University, \\
\hphantom{$^\ddag$}~1-13-27 Kasuga, Bunkyo-ku Tokyo 112-8551, Japan}
\EmailD{\href{mailto:takemura@math.chuo-u.ac.jp}{takemura@math.chuo-u.ac.jp}}

\ArticleDates{Received July 24, 2012, in f\/inal form November 07, 2012; Published online November 15, 2012}

\Abstract{Inf\/initely many explicit solutions of certain second-order dif\/ferential equations with an apparent singularity of
 characteristic exponent $-2$ are constructed by adjusting the parameter of the multi-indexed
Laguerre polynomials.}


\Classification{33C45; 33C47}

\renewcommand{\thefootnote}{\arabic{footnote}}
\setcounter{footnote}{0}

\section{Introduction}
\label{sec:intro}

In recent years, global solutions of Fuchsian dif\/ferential equations with more than three re\-gu\-lar singularities and those of the dif\/ferential equations of conf\/luent types, attracted attention~\cite{oshima2, oshima}.
The  multi-indexed Laguerre and Jacobi polynomials~\cite{gomez3, os25}\footnote{For the bibliography of various exceptional orthogonal polynomials see  references in~\cite{os25}.}
and their simplest (one-indexed) cases, i.e.\ the exceptional orthogonal polynomials,
provide explicit examples of inf\/initely many global solutions of such equations.
They are obtained by Darboux--Crum transformations \cite{adler,crum, darb,krein}
applied to the exactly solvable Schr\"odinger equations \cite{susyqm, infhul}
of the radial oscillator potential   and the P\"oschl--Teller potential.
In these global solutions, all the extra singularities
(i.e.\ those added by Darboux--Crum transformations) turn out to be apparent and the characteristic
exponents are all $-1$. This is a generic property of the Darboux--Crum transformations.
The role of Darboux transformations in iso-monodromic deformations has been known
for some time \cite{dui-gr}  and some explicit examples of monodromy free rational potentials were
given in \cite{gibbons, oblom}.

In the present paper,  we construct, by adjusting the parameter of the multi-indexed Laguerre polynomials,
several families of inf\/initely many global solutions  of Schr\"odinger equations with two or more
f\/inite re\-gu\-lar singularities and one irregular singularity at inf\/inity.
At one of the re\-gu\-lar singular points, these solutions have an apparent singularity of characteristic exponent~$-2$.
We now show one of the examples which is labeled as~(E) in Section~\ref{sec:milag}.
The Hamiltonian is given by
\begin{equation*}
\mathcal{H}=-\frac{d^2}{dx^2}+x^2+\frac{195}{4x^2}-16 -2\frac{d^2\log (x^2+6 )^3(x^2+14)}{dx^2}.
\end{equation*}
We change the variable by $\eta =x^2$.
Then the equation $\mathcal{H}\psi(x)=\mathcal{E}\psi(x)$ $(\mathcal{E} \in\mathbb{C})$ in the va\-riab\-le~$\eta $ has apparent singularities $\eta = -6$ and $\eta = -14$, whose characteristic exponent is~$-2$ and~$-1$ respectively.
All the apparent singularities are located outside of the half line $(0,\infty)$, thus the eigenfunctions of the Schr\"odinger equations are constructed by using Darboux--Crum transformation, and they are global solutions of the dif\/ferential equations.
The essential part of the numerators of these global solutions form orthogonal polynomials over a half line $(0,\infty)$, with the square of the denominator forming the part of the weight functions $e^{-x}x^{7}/\{ (x+6)^{4}(x+14)^{2} \}$.
The polynomials in the numerators of the global solutions are obtained within the framework of multi-indexed orthogonal polynomials, and their degrees have ``gaps''.
In this example, the polynomials of degrees~0,~1, and~2 are missing.
We also report examples (A), (B), (C), (D), (F), (G) in the two index case and an example in the three index case.
Each example has its peculiar properties (see section three).
The orthogonal polynomials related to (A) and (B) are also reported in~\cite{gomez4}.

This paper is organised as follows. In Section~\ref{sec:darb}, we review the Darboux transformations with the induced singularities.
In Section~\ref{sec:milag} we show that the special cases of multi-indexed
Laguerre polynomials provide several families of inf\/initely many global solutions of conf\/luent
 dif\/ferential equations with two or more f\/inite regular singularities and one irregular singularity at
inf\/inity.
 They have an apparent singularity of characteristic exponent~$-2$. The f\/inal section is for a~summary and comments.

\section{Darboux transformations and induced singularities}
\label{sec:darb}

Our starting point is a generic Schr\"odinger
equation with a Hamiltonian $\mathcal{H}$
\begin{equation*}
\mathcal{H}=-\frac{d^2}{dx^2}+U(x),
\end{equation*}
in which $U(x)$ is a meromorphic potential.
Let $\psi(x)$ and $\varphi(x)$ be two distinct solutions,
not necessarily square-integrable eigenfunctions,  of
the Schr\"odinger equation
\begin{equation}
\mathcal{H}\psi(x)=\mathcal{E}\psi(x),
\qquad \mathcal{H}\varphi(x)= \tilde{\mathcal{E}}\varphi(x),\qquad
\mathcal{E}, \tilde{\mathcal{E}}\in\mathbb{C}.
\label{scheq}
\end{equation}
Then it is elementary to show that a new function
\begin{equation*}
{\psi}^{(1)}(x)\eqdef \frac{\text{W}[\varphi,\psi](x)}{\varphi(x)} = \frac{\varphi (x) \psi '(x) -\varphi '(x) \psi (x)}{\varphi(x)}
\end{equation*}
is a solution of a new Schr\"odinger equation with a deformed Hamiltonian ${\mathcal H}^{(1)}$
\begin{equation*}
{\mathcal H}^{(1)}=-\frac{d^2}{dx^2}+{U}^{(1)}(x),\qquad
{U}^{(1)}(x)\eqdef U(x)-2\frac{d^2\log\varphi(x)}{dx^2},
\end{equation*}
with the same energy $\mathcal{E}$
\begin{equation}
{\mathcal H}^{(1)}{\psi}^{(1)}(x)=\mathcal{E}{\psi}^{(1)}(x).
\label{newschr}
\end{equation}
The zeros of the ``seed'' solution $\varphi(x)$ produce singularities of the equation, and there are no multiplicities of the zeros of $\varphi(x)$ except for the singularities of $U(x)$, because the function $\varphi(x)$ satisf\/ies the Schr\"odinger equation and the exponents at a holomorphic point are~$0$ or~$1$.

To obtain the deformed Hamiltonian which admits an apparent singularity of exponent $-2$ or higher,
we use
Darboux--Crum transformations \cite{adler,crum, darb,krein}.
Let $\{\varphi_j(x), \tilde{\mathcal{E}}_j\}$, $j=1,\ldots,M$
be distinct seed solutions of the original Schr\"odinger equation~\eqref{scheq} and as above~$\psi(x) $ be a solution with the energy~$\mathcal{E}$
\begin{equation*}
\mathcal{H}\psi(x)=\mathcal{E}\psi(x),\qquad \mathcal{H}\varphi_j(x)
= \tilde{\mathcal{E}}_j\varphi_j(x),
\qquad \mathcal{E}, \tilde{\mathcal{E}}_j\in\mathbb{C},
\quad j=1,\ldots,M.
\end{equation*}
By repeating the above Darboux transformation $M$-times, we obtain a new function
\begin{equation}
\psi^{(M)}(x)\eqdef\frac{\text{W}[\varphi_1,\ldots,\varphi_M,\psi](x)}
{\text{W}[\varphi_1,\ldots,\varphi_M](x)},
\label{psiM}
\end{equation}
which satisf\/ies an $M$-th deformed Schr\"odinger equation
with the same energy \cite{adler, crum,krein}
\begin{gather}
{\mathcal H}^{(M)}{\psi}^{(M)}(x)=\mathcal{E}{\psi}^{(M)}(x),\nonumber\\
{\mathcal H}^{(M)}=-\frac{d^2}{dx^2}+{U}^{(M)}(x),\qquad
{U}^{(M)}(x)\eqdef U(x)-2\frac{d^2\log\text{W}[\varphi_1,\ldots,\varphi_M](x)}{dx^2},\nonumber\\
\big(\psi^{(M)},\psi^{(M)}\big)=\prod_{j=1}^M(\mathcal{E}- \tilde{\mathcal{E}}_j)\cdot(\psi,\psi).
\label{intformM}
\end{gather}
Here $\text{W}[f_1,\ldots,f_N](x)$ is a Wronskian
\begin{equation*}
\text{W}[f_1,\ldots,f_N](x)\eqdef
\text{Det}\left(\frac{d^{k-1}f_j(x)}{dx^{k-1}}\right)_{1\le j,k\le N},
\end{equation*}
and the inner product among real functions is def\/ined as usual
\begin{equation*}
(f,g)\eqdef\int  f(x)g(x)dx.
\end{equation*}
Except for the singularities of the potential, a zero of the Wronskian $\text{W}[\varphi_1,\ldots,\varphi_M](x)$
\begin{equation*}
\text{W}[\varphi_1,\ldots,\varphi_M](x) =c_0(x-x_0)^m+O(x-x_0)^{m+1},\qquad m\in\mathbb{Z}_{>0},
\end{equation*}
corresponds to an apparent singular point of  the deformed Hamiltonian ${\mathcal H}^{(M)}$
\begin{equation*}
x\approx x_0,\qquad {U}^{(M)}(x)=U(x_0)
+\frac{2m}{(x-x_0)^2}+\text{regular terms}.
\end{equation*}
The characteristic exponents of the solutions of the new Schr\"odinger equation \eqref{newschr} at the apparent singularity $x_0$ are determined by
\begin{equation*}
\rho(\rho-1)-2m=0.
\end{equation*}
For certain values of the integer $m$, the monodromy of the added singularity becomes trivial
\begin{gather*}
m=1,\quad \rho=-1,\; 2,\\
m=3,\quad \rho=-2,\; 3,\\
m=6,\quad \rho=-3,\; 4,\quad \text{etc}.
\end{gather*}
When all the added singularities are of trivial monodromy,
the possibility of global solutions arises.
The inf\/inite families of the exceptional (or one-indexed) Laguerre and Jacobi polyno\-mials~\cite{os25} are the examples of the global solutions of second-order Fuchsian (and its conf\/luent limit) equations with many extra singularities with all $\rho=-1$.
This is the generic case, since any  solution
$\varphi(x)$ of the Schr\"odinger equation~\eqref{scheq}  with the radial oscillator potential~\eqref{pot} or the P\"oschl--Teller potential~\cite{os25} cannot have a multiple zero ($m\ge2$).

\section{Multi-indexed Laguerre  polynomials}
\label{sec:milag}

In this section, we present several families of inf\/initely
many global solutions with a high deg\-ree ($\rho=-2$)
of apparent singularity. The strategy is simple.
We start from certain ``non-generic'' multi-indexed
Laguerre polynomials and choose the parameter $g$ such that
a  triple zero is achieved for the seed solution (the Wronskian).
Then we construct the global solutions expli\-citly.
When the zeros are outside the half line $(0,\infty)$, the global solutions are square integrable and they are eigenfunctions of the deformed Hamiltonian.

\subsection{``Non-generic" multi-indexed Laguerre polynomials}
\label{sec:nongen}

For the  multi-indexed
Laguerre polynomials \cite{os25}, we adopt the Schr\"odinger equation
\eqref{scheq} with the radial oscillator potential~\cite{susyqm, infhul}
\begin{equation}
U(x;g)=x^2+\frac{g(g-1)}{x^2}-(1+2g),\qquad g>0,\quad 0<x<\infty.
\label{pot}
\end{equation}
It has a regular singularity at $x=0$ with the characteristic exponents
$\rho=g,1-g$ and an irregular singularity at $x=\infty$.
The eigenpolynomial solutions are
\begin{gather}
\mathcal{H}\phi_n(x;g) =\mathcal{E}_n\phi_n(x;g),
\qquad   \mathcal{E}_n=4n,\nonumber\\
\phi_n(x;g)=e^{-x^2/2}x^gL_n^{\big(g-\tfrac12\big)}(\eta(x)),
\qquad \eta(x)\eqdef x^2,\qquad n=0,1,\ldots,
\label{lagint}
\end{gather}
in which $L_n^{(\alpha)}(\eta)$ is the degree $n$ Laguerre polynomial
in $\eta$.
They form a complete set of ortho\-go\-nal polynomials over $(0,\infty)$
\begin{gather*}
   (\phi_n,\phi_m) \eqdef \int_0^\infty\!  \phi_n(x;g)\phi_m(x;g)dx =
  \frac12 \int_0^\infty\!  e^{-\eta}\eta^{g-\tfrac12}L_n^{\big(g-\tfrac12\big)}(\eta)L_m^{\big(g-\tfrac12\big)}(\eta)d\eta
=h_n(g)\delta_{nm},\nonumber\\
    h_n(g)\eqdef\frac{1}{2n!}\Gamma\big(n+g+\tfrac12\big).
\end{gather*}
The generic multi-indexed Laguerre polynomials are
constructed by the 
Darboux--Crum transformations~\eqref{psiM}
in terms of  two types of polynomial
($\text{v}\in\mathbb{Z}_{>0}$) seed solutions, called the virtual state solutions~\cite{os25},
\begin{alignat}{4}
  &  \tilde{\phi}_\text{v}^\ai(x;g) \eqdef
  e^{\frac12x^2}x^gL_\text{v}^{\big(g-\tfrac12\big)}(-\eta(x)), \quad&&
  \tilde{\mathcal{E}}_\text{v}^\ai(g)\eqdef-4\big(g+\text{v}+\tfrac12\big),\quad &&
   \big(\tilde{\phi}_\text{v}^\ai, \tilde{\phi}_\text{v}^\ai\big)=\infty, \label{seed1} &
\\
&  \tilde{\phi}_\text{v}^\ait(x;g) \eqdef
  e^{-\frac12x^2}x^{1-g}L_\text{v}^{\big(\tfrac12-g\big)}(\eta(x)), \quad&&
   \tilde{\mathcal{E}}_\text{v}^\ait(g)\eqdef-4\big(g-\text{v}-\tfrac12\big),
   \quad && \big(\tilde{\phi}_\text{v}^\ait, \tilde{\phi}_\text{v}^\ait\big)=\infty, \label{seed2}  &
   \end{alignat}
with certain constraints on the range of the parameter~$g$.
In order to generate a wider class of global solutions
with a high degree of apparent singularities,
 we will adopt {\em non-virtual state} polynomial  seed solutions, too:
\begin{gather}
    \tilde{\phi}_\text{v}^\aitr(x;g) \eqdef
  e^{\frac12x^2}x^{1-g}L_\text{v}^{\big(\tfrac12-g\big)}(-\eta(x)), \qquad
  \tilde{\mathcal{E}}_\text{v}^\aitr(g)\eqdef-4(\text{v}+1),
  \qquad \big(\tilde{\phi}_\text{v}^\aitr, \tilde{\phi}_\text{v}^\aitr\big)=\infty.\!\!\!
   \label{seed3}
\end{gather}
Here we summarise the def\/inition of ``non-generic'' multi-indexed Laguerre polynomials.

\begin{definition*}
The  ``non-generic'' multi-indexed Laguerre polynomials are
the main part of the eigenfunctions $\psi^{(M)}(x)$ \eqref{psiM} generated
by the multiple Darboux transformations with the eigenfunctions $\{\phi_n(x;g)\}$
and the three types of seed solutions \eqref{seed1}, \eqref{seed2}
and \eqref{seed3}, the virtual and non-virtual state solutions,
and that the parameter $g$ is restricted only by the positivity $g>0$.
\end{definition*}

The regularity of the solutions and/or the positivity of the resultant weight
functions etc. must be verif\/ied in each specif\/ic case of the chosen parameter values. This is in good contrast with the multi-indexed Laguerre and Jacobi polynomials def\/ined in \cite{os25}, in which the parameters are
restricted (equations~(23) and~(24) of~\cite{oshima2}) so that the seed solutions are of def\/inite sign and their inverse are not square integrable
\begin{equation*}
\big(1/\tilde{\phi}_\text{v}^\ai, 1/\tilde{\phi}_\text{v}^\ai\big),
=\big(1/\tilde{\phi}_\text{v}^\ait, 1/\tilde{\phi}_\text{v}^\ait\big)=\infty.
\end{equation*}
These two additional properties guarantee the regularity of the solutions and the positivity of the weight function. In other words,
these restrictions are lifted for the possible construction of the
non-generic multi-indexed Laguerre polynomials.
In the singular solution cases, the square integrability
is irrelevant and we also consider  $g<0$ example~(F).

For the $M$-index case, the seed solutions is  $\text{W}[\varphi_{\text{v}_1},\ldots,\varphi_{\text{v}_M}](x)$.
We look for the possibility that this function has  triple  or higher zeros at some points.

\subsection{Two-index case}
\label{sec:two-ind}

Here we consider the simplest multi-index case, the two-index case
\begin{equation}
\phi^{(2)}_n(x)\eqdef
\frac{\text{W}[\varphi_1,\varphi_2,\phi_n](x)}{\text{W}[\varphi_1,\varphi_2](x)}.
\label{phi2}
\end{equation}
We focus on the seed solutions of the form
\begin{equation*}
\text{W}\big[ \tilde{\phi}_{\text{v}_1}^\aitr(x;g),
 \tilde{\phi}_{\text{v}_2}^\ai(x;g)\big](x),\qquad
\text{W}\big[ \tilde{\phi}_{\text{v}_1}^\aitr(x;g),
\tilde{\phi}_{\text{v}_2}^\ait(x;g)\big](x),\qquad
\text{W}\big[ \tilde{\phi}_{\text{v}_1}^\ai(x;g),
\tilde{\phi}_{\text{v}_2}^\ait(x;g)\big](x),
\end{equation*}
with lower degrees $\text{v}_1$ and $\text{v}_2$
and evaluate the discriminants of the polynomial
part of the seed solutions.   At the roots of the discriminant, which is a polynomial in~$g$,
the above seed solution has a multiple zero.
The f\/irst two are the non-generic type and the last one is the generic type. As we will show in some detail,
they have dif\/ferent features.

So far we have encountered the following seven cases:
\begin{alignat*}{3}
&(A): \ \text{W}\big[ \tilde{\phi}_{1}^\aitr\big(x;\tfrac34\big),
\tilde{\phi}_{2}^\ai\big(x;\tfrac34\big)\big](x),\qquad && (B): \ \text{W}\big[ \tilde{\phi}_{2}^\aitr\big(x;\tfrac14\big),
\tilde{\phi}_{1}^\ai\big(x;\tfrac14\big)\big](x),& \\
&(C):\ \text{W}\big[ \tilde{\phi}_{2}^\aitr\big(x;\tfrac94\big),
\tilde{\phi}_{1}^\ait\big(x;\tfrac94\big)\big](x), \qquad &&
 (D):\ \text{W}\big[ \tilde{\phi}_{1}^\aitr\big(x;\tfrac94\big),
\tilde{\phi}_{2}^\ait\big(x;\tfrac94\big)\big](x), & \\
&(E):\ \text{W}\big[ \tilde{\phi}_{2}^\ai\big(x;\tfrac{15}{2}\big),
\tilde{\phi}_{1}^\ait\big(x;\tfrac{15}{2}\big)\big](x),\qquad &&
 (F):\ \text{W}\big[ \tilde{\phi}_{1}^\ai\big(x;-\tfrac{13}{2}\big),
\tilde{\phi}_{2}^\ait\big(x;-\tfrac{13}{2}\big)\big](x),& \\
&(G):\ \text{W}\big[ \tilde{\phi}_{3}^\ai\big(x;\tfrac{53}{2}\big),
\tilde{\phi}_{1}^\ait\big(x;\tfrac{53}{2}\big)\big](x).\qquad &&&
\end{alignat*}

\subsubsection{Seed solution case~(A)}
\label{sec:seedtwo-A}

 In this case
\begin{gather*}
 \text{W}\big[ \tilde{\phi}_{1}^\aitr(x;g), \tilde{\phi}_{2}^\ai(x;g)\big](x) \\
 \qquad{}
= \tfrac{1}{16} (2 g+1)  e^{x^2}
\big({-}9+18 g+4 g^2-8 g^3+\big(18-8 g^2\big) x^2+(12+8 g) x^4+8 x^6\big).
\end{gather*}
The discriminant of the polynomial in the variable $\eta (=x^2)$ is
\begin{equation*}
2048 (2 g-3) (2 g+3)^2 (4 g-3)^2,
\end{equation*}
and $g=\tfrac34$ gives a cubic zero
\begin{equation*}
 \text{W}\big[ \tilde{\phi}_{1}^\aitr\big(x;\tfrac34\big),
 \tilde{\phi}_{2}^\ai\big(x;\tfrac34\big)\big](x)=
 \tfrac{5}{256} e^{x^2}\big(3+4 x^2\big)^3.
\end{equation*}
The potential is
\begin{equation}
U_A(x)=x^2-\frac{3}{16 x^2}+\frac{48}{3+4 x^2}
-\frac{288}{\big(3+4 x^2\big)^2}-\frac{13}{2}.
\label{potA}
\end{equation}

\subsubsection{Seed solution case (B)}
\label{sec:seedtwo-B}

 In this case
 \begin{gather*}
\text{W}\big[ \tilde{\phi}_{2}^\aitr(x;g), \tilde{\phi}_{1}^\ai(x;g)\big](x)\\
\qquad{}=\tfrac{1}{16} (2 g-3)  e^{x^2}   \big(5-2 g-20g^2+8 g^3+\big(10+16g-8g^2\big)x^2
  +(20 -8g) x^4+8 x^6\big).
\end{gather*}
The discriminant of the polynomial in $\eta $ is
\begin{equation*}
-2048 (2 g-5)^2 (2 g+1) (4 g-1)^2,
\end{equation*}
and $g=\tfrac14$ gives a cubic zero
\begin{equation*}
 \text{W}\big[ \tilde{\phi}_{2}^\aitr\big(x;\tfrac14\big), \tilde{\phi}_{1}^\ai\big(x;\tfrac14\big)\big](x)=
-\tfrac{5}{256} e^{x^2} \big(3+4 x^2\big)^3.
\end{equation*}
The potential is the same as $U_A(x)$ up to an additive constant
\begin{equation}
U_B(x)=U_A(x)+1=x^2-\frac{3}{16 x^2}+\frac{48}{3+4 x^2}
-\frac{288}{\big(3+4 x^2\big)^2}-\frac{11}{2}.
\label{potB}
\end{equation}

\subsubsection{Seed solution case (C)}
\label{sec:seedtwo-C}
 In this case
 \begin{gather}
 \text{W}\big[ \tilde{\phi}_{2}^\aitr(x;g), \tilde{\phi}_{1}^\ait(x;g)\big](x)\nonumber\\
 {}=\tfrac{1}{8} x^{3-2 g} \big({-}135+174 g-68 g^2+8 g^3
 +\big({-}54+48 g -8 g^2\big) x^2+(36-8 g) x^4 +8 x^6\big).
 \label{wC}
\end{gather}
The discriminant of the polynomial in $\eta $ is
\begin{equation*}
-2048 (2 g-9)^2 (2 g-3) (4 g-9)^2,
\end{equation*}
and $g=\tfrac94$ gives a cubic zero
\begin{equation*}
\text{W}\big[ \tilde{\phi}_{2}^\aitr\big(x;\tfrac94\big),
\tilde{\phi}_{1}^\ait\big(x;\tfrac94\big)\big](x)=
\frac{\big(3 + 4 x^2\big)^3}{64 x^{3/2}}.
\end{equation*}
The potential is exactly the same as $U_B$
\begin{equation*}
U_C(x)=U_B(x)=x^2-\frac{3}{16 x^2}+\frac{48}{3+4 x^2}
-\frac{288}{\big(3+4 x^2\big)^2}-\frac{11}{2}.
\end{equation*}
Since the eigenvalue $\mathcal{E}_n=4n$ is independent of~$g$, the corresponding
solutions $\phi_n^{(2)}(x)$~\eqref{phi2} for~(B) and~(C) satisfy the same equations.
Thus they are {\em identical}. We will not consider  case~(C) in the next subsection.

\subsubsection{Seed solution case (D)}
\label{sec:seedtwo-D}

In this case
\begin{gather*}
\text{W}\big[\tilde{\phi}_{1}^\aitr(x;g), \tilde{\phi}_{2}^\ait(x;g)\big](x) \\
\qquad{}=\tfrac{1}{8} x^{3-2 g} \big({-}135+174 g-68 g^2+8 g^3
 +\big(54-48 g +8 g^2\big) x^2+(36-8 g) x^4 -8 x^6\big).
\end{gather*}
The polynomial part is obtained from that of \eqref{wC}
by the change $x^2\to -x^2$.
Thus the discri\-mi\-nant is the same and $g=\tfrac94$ gives a cubic zero
\begin{equation*}
\text{W}\big[ \tilde{\phi}_{1}^\aitr\big(x;\tfrac94\big), \tilde{\phi}_{2}^\ait\big(x;\tfrac94\big)\big](x)=-
\frac{\big({-}3 + 4 x^2\big)^3}{64 x^{3/2}}.
\end{equation*}
The potential has a singularity in $(0,\infty)$
\begin{equation}
U_D(x)=x^2-\frac{3}{16 x^2}+\frac{48}{-3+4 x^2}
+\frac{288}{\big({-}3+4 x^2\big)^2}-\frac{11}{2}.
\label{potD}
\end{equation}
The corresponding solutions \eqref{phi2} are not square integrable
and  this sequence is not that of  orthogonal polynomials.

Cases (E)--(G) belong to the generic type.

\subsubsection{Seed solution case (E)}
\label{sec:seedtwo-E}
In this case
\begin{gather*}
\text{W}\big[\tilde{\phi}_{2}^\ai(x;g), \tilde{\phi}_{1}^\ait(x;g)\big](x)\\
= -\tfrac{9}{16}+\tfrac{5}{2} g^2-g^4+\left(\tfrac{9}{2}+9g-2g^2-4g^3\right)x^2 +\left(\tfrac{15}{2}-4g-6g^2\right) x^4-2(1+2g) x^6-x^8.
\end{gather*}
The discriminant of the polynomial in $\eta $ is
\begin{gather*}
-16 (2 g-15)^2 (2 g -3) (2 g +1) (2 g +3)^2
\end{gather*}
 and $g=\tfrac{15}{2}$ gives a cubic zero
\begin{equation*}
\text{W}\big[ \tilde{\phi}_{2}^\ai\big(x;\tfrac{15}{2}\big),
\tilde{\phi}_{1}^\ait\big(x;\tfrac{15}{2}\big)\big](x)=
\big(6+ x^2\big)^3\big(14+x^2\big).
\end{equation*}
The potential is
\begin{equation}
U_E(x)= x^2+\frac{195}{4 x^2}-\frac{144}{\big( 6 + x^2 \big)^2}+\frac{12}{6+x^2}-\frac{112}{\big( 14+x^2 \big)^2}+\frac{4}{14 +x^2}-16.
\label{potE}
\end{equation}

\subsubsection{Seed solution case (F)}
\label{sec:seedtwo-F}
In this case
\begin{gather*}
\text{W}\big[\tilde{\phi}_{1}^\ai(x;g), \tilde{\phi}_{2}^\ait(x;g)\big](x)
=\tfrac{15}{16}-g -\tfrac{7}{2}g^2+4g^3-g^4 -\left(\tfrac{15}{2}+7g-14g^2+4g^3\right)x^2\\
\qquad {}-\left(\tfrac{5}{2}-16 g+6 g^2\right)x^4
+(6-4 g)x^6-x^8.
\end{gather*}
The discriminant of the polynomial in $\eta $ is
\begin{gather*}
-16 ( 2 g -5)^2 (2 g -3) (2 g +1) (2 g +13)^2
\end{gather*}
 and at a negative value of the parameter $g=-\tfrac{13}{2}$ it gives a cubic zero
\begin{equation*}
\text{W}\big[\tilde{\phi}_{1}^\ai\big(x;-\tfrac{13}{2}\big), \tilde{\phi}_{2}^\ait\big(x;-\tfrac{13}{2}\big)\big](x)=
-\big(x^2-14\big) \big(x^2-6\big)^3.
\end{equation*}
The potential is singular in $(0,\infty)$
\begin{equation}
U_F(x)=  x^2+\frac{195}{4 x^2}+\frac{144}{\big(x^2-6\big)^2}
+\frac{12}{x^2-6}+\frac{112}{\big(x^2-14\big)^2}
+\frac{4}{x^2-14} +12.
\label{potF}
\end{equation}

\subsubsection{Seed solution case (G)}
\label{sec:seedtwo-G}
In this case
\begin{gather*}
 \text{W}\big[\tilde{\phi}_{3}^\ai(x;g), \tilde{\phi}_{1}^\ait(x;g)\big](x)
=\tfrac{1}{96}\big\{
 (-3+2g)(-1+2g)(1+2g)(3+2g)(5+2g) \\
\qquad{}+10(-3+2g)(1+2g)(3+2g)(5+2g)x^2+80(-1+g)(3+2g)(5+2g)x^4\\
\qquad {} +160g(5+2g)x^6+80(3+2g)x^8+32x^{10}\big\}.
\end{gather*}
The discriminant in $\eta $ is proportional to
\begin{gather*}
 (2 g -3 ) ( 2 g +1) (2 g +3)^2 (2 g +5)^3 ( 10 g -39 )^2
\end{gather*}
and $g=\tfrac{39}{10}$ gives a cubic zero
\begin{gather*}
 \text{W}\big[\tilde{\phi}_{3}^\ai\big(x;\tfrac{39}{10}\big), \tilde{\phi}_{1}^\ait\big(x;\tfrac{39}{10}\big)\big](x)=\tfrac{1}{9375}\big(12+5x^2\big)^3\big(2244+495x^2+25x^4\big).
\end{gather*}
The potential is
\begin{gather*}
U_G(x) = x^2+\frac{1131}{100 x^2}-\frac{44}{5}
-\frac{1440}{\big(12+5
   x^2\big)^2}+\frac{60}{12+5 x^2}\\
\hphantom{U_G(x) =}{}
  +\frac{165000 x^2}{\big(2244+495 x^2+25
   x^4\big)^2}+\frac{20 \big(10 x^2-99\big)}{2244+495 x^2+25
   x^4} .
\end{gather*}
The basic structure of~(G)  is essentially the same as that of (E). So we will not present the explicit formulas of~(G), since they tend to be rather
lengthy.

\subsubsection{The detailed structure of these seed solutions}
\label{sec:seedtwo-str}

In the non-generic cases (A) and (B), the potential functions
are essentially the same, but the main
eigenfunctions carry dif\/ferent coupling constants,
$\big\{\phi_n\big(x;\tfrac34\big)\big\}$ and $\big\{\phi_n\big(x;\tfrac14\big)\big\}$.
Thus the corresponding multi-indexed polynomials are entirely dif\/ferent.
The equality (up to a sign) of the seed solutions of
(A) and (B) is easily understood,
when one notes the  trivial identity
\begin{equation*}
\tilde{\phi}_{\text v}^\ai(x;g)=\tilde{\phi}_{\text v}^\aitr(x;1-g).
\end{equation*}
In fact, we have
\begin{gather}
  \tilde{\phi}_{1}^\ai\big(x;\tfrac14\big)= \tilde{\phi}_{1}^\aitr\big(x;\tfrac34\big)=\tfrac14e^{x^2/2}x^{1/4}\big(3+4x^2\big),\nonumber \\
  \big(1/\tilde{\phi}_{1}^\ai\big(x;\tfrac14\big),1/\tilde{\phi}_{1}^\ai\big(x;\tfrac14\big)\big)
  =\big(1/\tilde{\phi}_{1}^\aitr\big(x;\tfrac34\big),1/\tilde{\phi}_{1}^\aitr\big(x;\tfrac34\big)\big)<\infty,
  \label{finnorm}\\
\tilde{\phi}_{2}^\ai\big(x;\tfrac34\big)= \tilde{\phi}_{2}^\aitr\big(x;\tfrac14\big)=\tfrac{1}{32}e^{x^2/2}x^{3/4}\big(3+4x^2\big)\big(15+4x^2\big), \nonumber\\
 \big(1/\tilde{\phi}_{2}^\ai\big(x;\tfrac34\big),1/\tilde{\phi}_{2}^\ai\big(x;\tfrac34\big)\big)
 =\big(1/\tilde{\phi}_{2}^\aitr\big(x;\tfrac14\big),1/\tilde{\phi}_{1}^\aitr\big(x;\tfrac14\big)\big)
 =\infty.\nonumber
\end{gather}
A cubic zero of the seed solution in (A) and (B) is realised by a nice magic of the Wronskian
\begin{gather*}
\text{W}\big[e^{x^2/2}x^{1/4}\big(3+4x^2\big), e^{x^2/2}x^{3/4}\big(3+4x^2\big)\big(15+4x^2\big)\big](x)\\
\qquad{}=\big(e^{x^2/2}x^{1/4}\big(3+4x^2\big)\big)^2\text{W}\big[1,x^{1/2}\big(15+4x^2\big)\big](x).
\end{gather*}
The Wronskian on the r.h.s.\ provides another factor $3+4x^2$.
The fact that $1/\tilde{\phi}_{1}^\ai\big(x;\tfrac14\big)$ and
$1/\tilde{\phi}_{1}^\aitr\big(x;\tfrac34\big)$ are square integrable \eqref{finnorm}
means that  the mechanism of multi-index polynomials developed in~\cite{os25} does not apply.

For  (C) and (D), the functions are
\begin{gather*}
  \tilde{\phi}_{1}^\ait\big(x;\tfrac94\big) =-\tfrac14e^{-x^2/2}x^{-5/4}\big(3+4x^2\big),
  \qquad
\tilde{\phi}_{2}^\ait\big(x;\tfrac94\big)=\tfrac{1}{32}e^{-x^2/2}x^{-5/4}\big({-}3+4x^2\big)\big(1+4x^2\big), \nonumber\\
   \tilde{\phi}_{1}^\aitr\big(x;\tfrac94\big) =\tfrac14e^{x^2/2}x^{-5/4}\big({-}3+4x^2\big),
   \qquad   \tilde{\phi}_{2}^\aitr\big(x;\tfrac94\big)
   =\tfrac{1}{32}e^{x^2/2}x^{-5/4}\big(3+4x^2\big)\big({-}1+4x^2\big), \\
 \big(\tilde{\phi}_{1}^\ait\big(x;\tfrac94\big),\tilde{\phi}_{1}^\ait\big(x;\tfrac94\big)\big)=\infty,
 \qquad \big(1/\tilde{\phi}_{1}^\ait\big(x;\tfrac94\big),1/\tilde{\phi}_{1}^\ait\big(x;\tfrac94\big)\big)=\infty.\nonumber
\end{gather*}
Again the  cubic zeros of these cases are realised by a nice magic of the Wronskian.
As expected, $\tilde{\phi}_{2}^\ait\big(x;\tfrac94\big)$,
$\tilde{\phi}_{1}^\aitr\big(x;\tfrac94\big)$ and
$\tilde{\phi}_{2}^\aitr\big(x;\tfrac94\big)$ change sign in $(0,\infty)$,
which is a deviation from the mechanism of multi-index polynomials
developed in~\cite{os25}.

Now it is clear that the solutions $\phi^{(2)}_n(x)$ \eqref{phi2}
have a $\rho=-2$ apparent singularity, that is,
\begin{equation}
\phi^{(2)}_n(x)\propto \frac{1}{(3+4x^2)^2},\quad \text{(A)},\ \text{(B)},\ \text{(C)},\qquad
\phi^{(2)}_n(x)\propto \frac{1}{(-3+4x^2)^2},\quad \text{(D)}.
\label{rho-2}
\end{equation}
Among the six terms of the numerator Wronskian
$\text{W}[\varphi_1,\varphi_2,\phi_n](x)$ of \eqref{phi2},
the four terms containing
$\phi_n'(x)$ or $\phi_n''(x)$ contain at least one term $\varphi_1(x)$ or
$\varphi_2(x)$ without a derivative.
Thus they contain at least one factor of $3+4x^2$ in (A), (B), (C)
and $-3+4x^2$ in (D). The remaining two terms which contain
$\phi_n(x)$ are proportional to
\begin{equation*}
\varphi_1'(x)\varphi_2''(x)-\varphi_1''(x)\varphi_2'(x).
\end{equation*}
Since $\varphi_j(x)$ is a solution of the original Schr\"odinger
equation \eqref{scheq}, they satisfy
\begin{equation*}
\varphi_j''(x)=\big(U(x)-\tilde{\mathcal E}_j\big)\varphi_j(x),
\qquad j=1,2.
\end{equation*}
Thus the above terms also contain one factor of
$\varphi_1(x)$ or $\varphi_2(x)$, that is, $3+4x^2$
in (A), (B), (C) and $-3+4x^2$ in (D).
Thus one factor of the cubic zero of the seed function $\text{W}[\varphi_1,\varphi_2](x)$ is cancelled and we arrive at \eqref{rho-2}.

The relationship between the generic cases (E) and (F) has similarities with that between the non-generic cases (C)
and (D).
The functions in (E) are
\begin{gather*}
\tilde{\phi}^\ai_2\big(x;\tfrac{15}{2}\big) =\tfrac12e^{x^2/2}x^{15/2}\big(6+x^2\big)\big(12+x^2\big),\qquad
\tilde{\phi}^\ait_1\big(x;\tfrac{15}{2}\big)=-e^{-x^2/2}x^{-13/2}\big(6+x^2\big).
\end{gather*}
The functions in (F) are
\begin{gather*}
\tilde{\phi}^\ai_1\big(x;-\tfrac{13}{2}\big) =e^{x^2/2}x^{-13/2}\big({-}6+x^2\big),\qquad\!
\tilde{\phi}^\ait_1\big(x;-\tfrac{13}{2}\big)=\tfrac12e^{-x^2/2}x^{15/2}\big({-}6+x^2\big)\big({-}12+x^2\big).
\end{gather*}
Thus it is also clear that the solutions $\phi^{(2)}_n(x)$ \eqref{phi2}
have a $\rho=-2$ apparent singularity, that is,
\begin{equation*}
\phi^{(2)}_n(x)\propto \frac{1}{(6+x^2)^2},\quad \text{(E)},\qquad
\phi^{(2)}_n(x)\propto \frac{1}{(x^2-6)^2},\quad \text{(F)}.
\end{equation*}

\subsection[Global solutions with a $\rho=-2$ apparent singularity]{Global solutions with a $\boldsymbol{\rho=-2}$ apparent singularity}
\label{sec:global}

Here we present the inf\/initely many global solutions
with a $\rho=-2$ apparent singularity in explicit forms.
They are simply obtained by evaluating $\phi^{(2)}_n(x)$
\eqref{phi2} in each case (A)--(F),
to be denoted by $\phi^A_n(x)$, $\phi^B_n(x)$, \ldots, $\phi^F_n(x)$, respectively.
Most formulas look simpler as  functions of
$\eta\equiv \eta(x)=x^2$, since the eigenpolynomials of the
radial oscillator potential are the Laguerre polynomials
in $\eta$, \eqref{lagint}.

\subsubsection{Global solutions  (A)}
\label{sec:sol-A}

The global solutions are
\begin{gather}
\mathcal{H}^A\phi^A_n(x) =\mathcal{E}_n\phi^A_n(x),
\qquad \mathcal{H}^A\eqdef-\frac{d^2}{dx^2}+U_A(x),\qquad \mathcal{E}_n=4n,
\label{schA}\\
\phi^A_n(x) =\frac{e^{-\eta/2}\eta^{3/8}}{(3+4\eta)^2}\mathcal{L}^A_n(\eta),\qquad n=0,1,\ldots,
\label{Awavfun}
\end{gather}
in which  $\mathcal{L}^A_n(\eta)$ ($n\ge0$) is a degree $n+3$  polynomial in $\eta$ def\/ined by
\begin{gather*}
 \mathcal{L}^A_n(\eta)  \eqdef
 \big({-}117+156\eta+208\eta^2+64\eta^3\big)L_n^{(1/4)}(\eta)\\
\hphantom{\mathcal{L}^A_n(\eta)  \eqdef}{}
 -4n(3+4\eta)^2L_n^{(1/4)}(\eta)-4\eta(3+4\eta)(15+4\eta)
 \partial_\eta L_n^{(1/4)}(\eta).
\end{gather*}
In evaluating the Wronskian $\text{W}\big[ \tilde{\phi}_{1}^\aitr\big(x;\tfrac34\big), \tilde{\phi}_{2}^\ai\big(x;\tfrac34\big),\phi_n\big(x;\tfrac34\big)\big](x)$, the second derivative of the Laguerre polynomial $L_n^{(1/4)}(\eta)$ appears.
It is replaced by $L_n^{(1/4)}(\eta)$ and $\partial_\eta L_n^{(1/4)}(\eta)$ by using the equation for the
Laguerre polynomial
\begin{equation*}
  x\partial_x^2L_n^{(\alpha)}(x)
  +(\alpha+1-x)\partial_xL_n^{(\alpha)}(x)
  +nL_n^{(\alpha)}(x)=0.
\end{equation*}
There is another member of this family corresponding to the eigenvalue
$\mathcal{E}=-8$, because of the fact that the Darboux transformation in terms
of $\tilde{\phi}_{1}^\aitr\big(x;\tfrac34\big)$ is not one-to-one,
$\big(1/\tilde{\phi}_{1}^\aitr\big(x;\tfrac34\big)$, $1/\tilde{\phi}_{1}^\aitr\big(x;\tfrac34\big)\big)<\infty$,
  \eqref{finnorm}, see~\cite{os25}.
Let us denote it tentatively by
$\mathcal{L}^A_{-2}(\eta)$
\begin{equation}
\mathcal{L}^A_{-2}(\eta)\eqdef 15+4\eta,
\label{agr}
\end{equation}
whose normalisation is irrelevant.
Note that $n=-1$ is missing.
  Some lower members are
\begin{gather*}
 \mathcal{L}^A_0(\eta) =-117+156\eta+208\eta^2+64\eta^3,\qquad
\mathcal{L}^A_1(\eta)= -\frac{765}{4} + 408 \eta + 408 \eta^2 - 64 \eta^4,\\
\mathcal{L}^A_2(\eta) =-\frac{8505}{32}+\frac{6237}{8} \eta
+567 \eta^2 -252 \eta^3-168 \eta^4+32 \eta^5.
\end{gather*}
They are also discussed in \S~6.2.5 of~\cite{gomez4}.
They satisfy the orthogonality relation
\begin{gather*}
 \int_0^\infty \frac{e^{-\eta}\eta^{\tfrac14}}{(3+4\eta)^4}\mathcal{L}^A_n(\eta) \mathcal{L}^A_m(\eta)d\eta  = h^A_n\delta_{nm},\qquad n,m=-2,0,1,\ldots, \\
 h^A_n \eqdef \frac{4(n+2)(4n+13)\Gamma(n+\frac54)}{n!},\qquad n\neq-2,
\end{gather*}
which is obtained by rewriting \eqref{intformM} with $\tilde{\mathcal E}^\aitr_1(\tfrac34)=-8$ and $\tilde{\mathcal E}^\ai_2(\tfrac34)=-13$.

By rewriting the Schr\"odinger equation for $\phi^A_n(x)$ \eqref{schA} in terms of
\begin{equation*}
\phi^A_n(x)=e^{-\eta/2}\eta^{3/8}y(\eta),
\end{equation*}
we obtain a second-order dif\/ferential equation with regular a singularity at $\eta=0$ and
$\eta=-\frac34$ and an irregular singularity at $x=\infty$
\begin{gather}
\eta y''+\left(\frac54-\eta\right)y'+\frac{45-24\eta+16\eta^2}{(3+4\eta)^2}y+n  y=0,\qquad n=-2,0,1, \ldots.
\label{difeqA}
\end{gather}
For each $n$ we have a global solution with a $\rho=-2$ apparent singularity
\begin{equation*}
y_n(\eta)=\frac{\mathcal{L}^A_n(\eta)}{(3+4\eta)^2},\qquad n=-2,0,1,\ldots .
\end{equation*}
The dif\/ferential equations for the polynomials $\{\mathcal{L}^A_n(\eta)\}$ read
\begin{gather*}
 4\eta(3+4\eta)\mathcal{L}^A_n(\eta)''-(4\eta-1)(15+4\eta)\mathcal{L}^A_n(\eta)'\n
\qquad {} +4\bigl(5+12\eta+n(3+4\eta)\bigr)\mathcal{L}^A_n(\eta)=0,\qquad
n=-2,0,1,\ldots .
\end{gather*}

The groundstate wavefunction $\phi^A_{-2}(x)$ \eqref{Awavfun} with $\mathcal{L}^A_{-2}$ \eqref{agr}
and the energy $\mathcal{E}=-8$
is related to the potential $U_A(x)$ \eqref{potA} through the prepotential
$w_A(x)$ (Riccati form of Schr\"odinger equation)
\begin{gather*}
e^{w_A(x)}\eqdef \frac{e^{-x^2/2}x^{3/4}(15+4x^2)}{(3+4x^2)^2}\quad
\Rightarrow \quad \left(\frac{dw_A(x)}{dx}\right)^2
+\frac{d^2w_A(x)}{dx^2}-8=U_A(x).
\end{gather*}
In other words, the Hamiltonian $\mathcal{H}^A$ \eqref{schA} is written in
a factorised form
\begin{gather}
\mathcal{H}^A =\mathcal{A}^{A\dagger}\mathcal{A}^{A}-8,\qquad
\mathcal{A}^{A}\eqdef \frac{d}{dx}-\frac{dw_A(x)}{dx},\qquad
\mathcal{A}^{A\dagger}=-\frac{d}{dx}-\frac{dw_A(x)}{dx}.
\label{HAsystem}
\end{gather}
It is easy to verify that the standard Darboux--Crum transformation
\cite{crum}
of deleting the groundstate applied to
\eqref{HAsystem}  will lead to the system of one-indexed
(or the exceptional) polynomial system generated in terms of the virtual
state solution $\tilde{\phi}^\ai_2\big(x;\tfrac34\big)$
\begin{gather*}
\mathcal{H}^{A1} \eqdef\mathcal{A}^{A}\mathcal{A}^{A\dagger}-8
=-\frac{d^2}{dx^2}+\left(\frac{dw_A(x)}{dx}\right)^2
-\frac{d^2w_A(x)}{dx^2}-8 =-\frac{d^2}{dx^2}+U_{A1}(x),\\
U_{A1}(x)
=
 x^2+\frac{21}{16x^2}+\frac{16}{3+4 x^2}
 -\frac{96}{\left(3+4x^2\right)^2}+\frac{16}{15+4 x^2}
   -\frac{480}{\left(15+4x^2\right)^2}
   -\frac{9}{2} \nonumber \\
\hphantom{U_{A1}(x)}{}
= U\big(x;\tfrac34\big)-2\frac{d^2}{dx^2}\log\tilde{\phi}^{\ai}_2\big(x;\tfrac34\big).
\end{gather*}
This shows clearly that the non-generic two-indexed polynomial system (A)
generated by a non-virtual state solution $\tilde{\phi}^\aitr_1\big(x;\tfrac34\big)$
and a virtual state solution $\tilde{\phi}^\ai_2\big(x;\tfrac34\big)$ is not
{\em shape-invariant} \cite{genden}, in contrast with the generic
multi-indexed polynomial systems \cite{os25}.
The standard Darboux--Crum transformation removes the groundstate
``created" by the non-virtual state solution $\tilde{\phi}^\aitr_1\big(x;\tfrac34\big)$
and the generic one-indexed polynomial system generated by the
virtual state solution $\tilde{\phi}^\ai_2\big(x;\tfrac34\big)$ is obtained.

\subsubsection{Global solutions  (B)}\label{sec:sol-B}

The global solutions are
\begin{gather}
\mathcal{H}^B\phi^B_n(x) =\mathcal{E}_n\phi^B_n(x),\qquad \mathcal{H}^B\eqdef-\frac{d^2}{dx^2}+U_B(x),\qquad \mathcal{E}_n=4n,
\label{schB}\\
\phi^B_n(x) =\frac{e^{-\eta/2}\eta^{1/8}}{(3+4\eta)^2}\mathcal{L}^B_n(\eta),\qquad n=0,1,\ldots,
\label{Bwavfun}
\end{gather}
in which  $\mathcal{L}^B_n(\eta)$ ($n\ge0$) is a degree $n+3$  polynomial
in $\eta$ def\/ined by
\begin{gather*}
 \mathcal{L}^B_n(\eta)  \eqdef
 \big({-}63+252\eta+240\eta^2+64\eta^3\big)L_n^{(-1/4)}(\eta)\\
\hphantom{\mathcal{L}^B_n(\eta)  \eqdef}{}
-4n(3+4\eta)^2L_n^{(-1/4)}(\eta)-4\eta(3+4\eta)(15+4\eta)\partial_\eta L_n^{(-1/4)}(\eta).
\end{gather*}
There is another member of this family corresponding to the eigenvalue
$\mathcal{E}=-12$, because of the fact that the Darboux transformation
in terms of $\tilde{\phi}_{2}^\aitr\big(x;\tfrac14\big)$
is not one to one,
$\big(1/\tilde{\phi}_{2}^\aitr\big(x;\tfrac14\big)$, $1/\tilde{\phi}_{2}^\aitr\big(x;\tfrac14\big)\big)<\infty$,
  \eqref{finnorm}, see \cite{os25}.
Let us  denote it tentatively by $\mathcal{L}^B_{-3}(\eta)$
\begin{equation}
\mathcal{L}^B_{-3}(\eta)\eqdef 1.
\label{Bgr}
\end{equation}
Its normalisation is irrelevant. Note that $n=-2,-1$ are missing.
Some lower-degree members are
\begin{gather*}
 \mathcal{L}^B_0(\eta) =-63+252\eta+240\eta^2+64\eta^3,\qquad
\mathcal{L}^B_1(\eta)= -\frac{297}{4} + 396 \eta
+ 264 \eta^2 - 64 \eta^3- 64 \eta^4,\\
\mathcal{L}^B_2(\eta) =-\frac{2835}{32}+\frac{4725}{8} \eta
+225 \eta^2 -300 \eta^3-120 \eta^4+32 \eta^5.
\end{gather*}
They are also discussed in \S~6.2.6 of~\cite{gomez4}.
They satisfy the orthogonality relation
\begin{gather*}
 \int_0^\infty \frac{e^{-\eta}\eta^{-1/4}}{(3+4\eta)^4}
 \mathcal{L}^B_n(\eta) \mathcal{L}^B_m(\eta)d\eta  = h^B_n\delta_{n m},
 \qquad n,m=-3,0,1,\ldots ,\\
 h^B_n \eqdef \frac{4(n+3)(4n+7)\Gamma(n+\frac34)}{n!}, \qquad n\neq-3,
\end{gather*}
which is obtained by rewriting \eqref{intformM} with
$\tilde{\mathcal E}^\aitr_2\big(\tfrac34\big)=-12$
and $\tilde{\mathcal E}^\ai_1\big(\tfrac34\big)=-7$.

By rewriting the Schr\"odinger equation for $\phi^B_n(x)$ \eqref{schB}
in terms of
\begin{equation*}
\phi^B_n(x)=e^{-\eta/2}\eta^{1/8}y(\eta),
\end{equation*}
we obtain a second-order dif\/ferential equation
with a regular singularity at $\eta=0$ and
$\eta=-\frac34$ and an irregular singularity at $x=\infty$
\begin{gather}
\eta y''+\left(\frac34-\eta\right)y'+\frac{45-24\eta+16\eta^2}{(3+4\eta)^2}y
+n  y=0,\qquad  n=-3,0,1, \ldots .
\label{difeqB}
\end{gather}
For each $n$ we have a global solution with a $\rho=-2$ apparent singularity
\begin{equation*}
y_n(\eta)=\frac{\mathcal{L}^B_n(\eta)}{(3+4\eta)^2},\qquad n=-3,0,1,\ldots .
\end{equation*}
The dif\/ferential equations for the polynomials $\{\mathcal{L}^B_n(\eta)\}$ read
\begin{gather*}
4\eta(3+4\eta)\mathcal{L}^B_n(\eta)''-\big({-}9+64\eta+16\eta^2\big)
\mathcal{L}^B_n(\eta)'\\
\qquad  {} +4(3+n)(3+4\eta)\mathcal{L}^B_n(\eta)=0,\qquad
n=-3,0,1,\ldots .
\end{gather*}

The groundstate wavefunction $\phi^B_{-3}(x)$ \eqref{Bwavfun} with $\mathcal{L}^B_{-3}$ \eqref{Bgr}
and the energy $\mathcal{E}=-12$
is related to the potential $U_B(x)$ \eqref{potB} through the prepotential
$w_B(x)$ (Riccati form of Schr\"odinger equation)
\begin{gather*}
e^{w_B(x)}\eqdef \frac{e^{-x^2/2}x^{1/4}}{(3+4x^2)^2}\quad
\Rightarrow \quad \left(\frac{dw_B(x)}{dx}\right)^2
+\frac{d^2w_B(x)}{dx^2}-12=U_B(x).
\end{gather*}
In other words, the Hamiltonian $\mathcal{H}^B$ \eqref{schB} is written in
a factorised form
\begin{gather*}
\mathcal{H}^B =\mathcal{A}^{B\dagger}\mathcal{A}^{B}-12,\qquad
\mathcal{A}^{B}\eqdef \frac{d}{dx}-\frac{dw_B(x)}{dx},\qquad
\mathcal{A}^{B\dagger}=-\frac{d}{dx}-\frac{dw_B(x)}{dx}.
\end{gather*}
It is easy to verify that the standard Darboux--Crum transformation~\cite{crum}
of deleting the groundstate will lead to the system of one-indexed
(or the exceptional) polynomial system generated in terms of the virtual
state solution $\tilde{\phi}^\ai_1\big(x;\tfrac14\big)$
\begin{gather*}
\mathcal{H}^{B1} \eqdef\mathcal{A}^{B}\mathcal{A}^{B\dagger}-12
=-\frac{d^2}{dx^2}+\left(\frac{dw_B(x)}{dx}\right)^2
-\frac{d^2w_B(x)}{dx^2}-12  =-\frac{d^2}{dx^2}+U_{B1}(x),\\
U_{B1}(x)
=
 x^2+\frac{5}{16x^2}+\frac{16}{3+4 x^2}
 -\frac{96}{\left(3+4x^2\right)^2}
   -\frac72  = U\big(x;\tfrac14\big)-2\frac{d^2}{dx^2}\log\tilde{\phi}^{\ai}_1\big(x;\tfrac14\big).
\end{gather*}

\subsubsection{Global solutions  (D)}\label{sec:sol-D}

The global solutions are
\begin{gather}
\mathcal{H}^D\phi^D_n(x) =\mathcal{E}_n\phi^D_n(x),
\qquad \mathcal{H}^D\eqdef-\frac{d^2}{dx^2}+U_D(x),
\qquad \mathcal{E}_n=4n,
\label{schD}\\
\phi^D_n(x)=\frac{e^{-\eta/2}\eta^{1/8}}{4(-3+4\eta)^2}
\mathcal{L}^D_n(\eta),\qquad n=0,1,\ldots,\nonumber
\end{gather}
in which  $\mathcal{L}^D_n(\eta)$ ($n\ge0$) is a degree $n+3$  polynomial in $\eta$ def\/ined by
\begin{gather*}
 \mathcal{L}^D_n(\eta)  \eqdef (3+4\eta)\big(63+16\eta^2\big)L_n^{(7/4)}(\eta)\\
\hphantom{\mathcal{L}^D_n(\eta)  \eqdef}{}
 -16n\eta(-3+4\eta)^2L_n^{(7/4)}(\eta)-36\eta(-3+4\eta)(1+4\eta)\partial_\eta L_n^{(7/4)}(\eta).
\end{gather*}
There is another member of this family belonging to the energy
$\mathcal{E}=-8$, to be denoted tentatively by $\mathcal{L}^D_{-2}(\eta)$
\begin{equation}
\mathcal{L}^D_{-2}(\eta)\eqdef 1+4\eta.
\label{Dgr}
\end{equation}
  Some lower members are
\begin{gather*}
 \mathcal{L}^D_0(\eta) =(3 + 4 \eta) \big(63 + 16 \eta^2\big),
 \qquad
\mathcal{L}^D_1(\eta)= \frac{2079}{4} + 792 \eta^2 - 384 \eta^3
 + 192 \eta^4,\\
\mathcal{L}^D_2(\eta) =\frac{31185}{32} - \frac{10395}{8} \eta
 + 3465 \eta^2 - 2940 \eta^3 + 1512 \eta^4 - 224 \eta^5.
\end{gather*}
Note that $n=-1$ is missing.
The lowest member $\mathcal{L}^D_{-2}$ \eqref{Dgr}
with the energy $\mathcal{E}=-8$
is related to the potential $U_D(x)$ \eqref{potD} through the prepotential
(Riccati form of Schr\"odinger equation)
\begin{gather*}
e^{w_D(x)}\eqdef \frac{e^{-x^2/2}x^{1/4}(1+4x^2)}{(-3+4x^2)^2}\quad
\Rightarrow \quad \left(\frac{dw_D(x)}{dx}\right)^2
+\frac{d^2w_D(x)}{dx^2}-8=U_D(x).
\end{gather*}

By rewriting the Schr\"odinger equation for $\phi^D_n(x)$ \eqref{schD}
in terms of
\begin{equation*}
\phi^D_n(x)=e^{-\eta/2}\eta^{1/8}y(\eta),
\end{equation*}
we obtain a second-order dif\/ferential equation
with a regular singularity at $\eta=0$ and
$\eta=+\frac34$ and an irregular singularity at $x=\infty$
\begin{gather}
\eta y''+\left(\frac34-\eta\right)y'-\frac{27+72\eta-16\eta^2}{(-3+4\eta)^2}y+n  y=0,\qquad   n=-2,0,1, \ldots .
\label{difeqD}
\end{gather}
For each $n$ we have a global solution with a $\rho=-2$ apparent singularity:
\begin{equation*}
y_n(\eta)=\frac{\mathcal{L}^D_n(\eta)}{(-3+4\eta)^2},\qquad n=-2,0,1,\ldots  .
\end{equation*}
The dif\/ferential equations for the polynomials $\{\mathcal{L}^D_n(\eta)\}$ read
\begin{gather*}
 4\eta(-3+4\eta)\mathcal{L}^D_n(\eta)''-(1+4\eta)(9+4\eta)
\mathcal{L}^D_n(\eta)'\\
\qquad {}  +4\bigl(3+12\eta+n(-3+4\eta)\bigr)\mathcal{L}^D_n(\eta)=0,\qquad
n=-2,0,1,\ldots  .
\end{gather*}

\subsubsection{Global solutions  (E)}\label{sec:sol-E}

The global solutions are
\begin{gather}
\mathcal{H}^E\phi^E_n(x) =\mathcal{E}_n\phi^E_n(x),
\qquad \mathcal{H}^E\eqdef-\frac{d^2}{dx^2}+U_E(x),\qquad \mathcal{E}_n=4n,
\label{schE}\\
\phi^E_n(x)=\frac{e^{-\eta/2}\eta^{\tfrac{15}{4}}}{(6+\eta)^2(14+\eta)}\mathcal{L}^E_n(\eta),\qquad n=0,1,\ldots,
\label{Ewavfun}
\end{gather}
in which  $\mathcal{L}^E_n(\eta)$ is a degree $n+3$  polynomial in $\eta$ def\/ined by
\begin{gather*}
 \mathcal{L}^E_n(\eta) \eqdef
 -24\big(840+280\eta+30\eta^2+\eta^3\big)L_n^{(7)}(\eta)\n
\hphantom{\mathcal{L}^E_n(\eta) \eqdef}{}  -4n(6+\eta)^2(14+\eta)L_n^{(7)}(\eta)-16\eta(6+\eta)(12+\eta)
 \partial_\eta L_n^{(7)}(\eta).
\end{gather*}
This belongs to the generic case and the Darboux transformation is one-to-one. The groundstate corresponds to $\mathcal{L}^E_0$.
  Some lower members are
\begin{gather}
 \mathcal{L}^E_0(\eta) =-24(840+280\eta+30\eta^2+\eta^3),\label{Egr}\\
\mathcal{L}^E_1(\eta) = 28\big({-}6336-1320\eta + 44 \eta^2 +22\eta^3+\eta^4\big),\nonumber\\
\mathcal{L}^E_2(\eta) =-16\big(54432+4536\eta-1944\eta^2-180\eta^3+12\eta^4+\eta^5\big).\nonumber
\end{gather}
They satisfy the orthogonality relation
\begin{gather*}
 \int_0^\infty \frac{e^{-\eta}\eta^{7}}{(6+\eta)^4(14+\eta)^2}\mathcal{L}^E_n(\eta) \mathcal{L}^E_m(\eta)d\eta = h^E_n\delta_{n m},\qquad n,m=0,1,\ldots, \\
 h^E_n\eqdef 16(n+10)(n+6)(n+1)_7,
\end{gather*}
which is obtained by rewriting \eqref{intformM} with $\tilde{\mathcal E}^\ai_2\big(\tfrac{15}{2}\big)=-40$ and $\tilde{\mathcal E}^\ait_1\big(\tfrac{15}{2}\big)=-24$. Here $(a)_n\eqdef\prod\limits_{k=1}^n(a+k-1)$ is the shifted factorial (the Pochhammer symbol).

By rewriting the Schr\"odinger equation for $\phi^E_n(x)$ \eqref{schE} in terms of
\begin{equation*}
\phi^E_n(x)=e^{-\eta/2}\eta^{15/4}y(\eta),
\end{equation*}
we obtain a second-order dif\/ferential equation with a regular singularity at $\eta=0$ and
$\eta=-6,-14$ and an irregular singularity at $x=\infty$
\begin{gather}
\eta y''+(8-\eta)y'+4\frac{1008+12\eta-16\eta^2-\eta^3}{(6+\eta)^2(14+\eta)^2}y+n y=0,\qquad   n=0,1, \ldots .
\label{difeqE}
\end{gather}
For each $n$ we have a global solution with a $\rho=-2$ and a $\rho=-1$ apparent singularities:
\begin{equation*}
y_n(\eta)=\frac{\mathcal{L}^E_n(\eta)}{(6+\eta)^2(14+\eta)},\qquad n=0,1,\ldots  .
\end{equation*}
The dif\/ferential equations for the polynomials $\{\mathcal{L}^E_n(\eta)\}$ read
\begin{gather*}
 \eta(6+\eta)(14+\eta)\mathcal{L}^E_n(\eta)''-\big({-}672-8\eta +18\eta^2 + \eta^3\big)\mathcal{L}^E_n(\eta)'\\
 \qquad  {}+\bigl({-}224+18\eta+3\eta^2+n(6+\eta)(14+\eta)\bigr)\mathcal{L}^E_n(\eta)=0,\qquad
n=0,1,\ldots  .
\end{gather*}

The groundstate wavefunction $\phi^E_{0}(x)$ \eqref{Ewavfun} with $\mathcal{L}^E_{0}$ \eqref{Egr}
is related to the potential $U_E(x)$ \eqref{potE} through the prepotential
$w_E(x)$:
\begin{gather*}
e^{w_E(x)}\eqdef \frac{e^{-x^2/2}x^{15/2}(840+280x^2+30x^4+x^6)}{(6+x^2)^2(14+x^2)}\quad
\Rightarrow \quad \left(\frac{dw_E(x)}{dx}\right)^2\!
+\frac{d^2w_E(x)}{dx^2}=U_E(x).
\end{gather*}
By construction, the system is shape-invariant
\begin{gather*}
 \left(\frac{dw_E(x)}{dx}\right)^2
-\frac{d^2w_E(x)}{dx^2} =U\big(x;\tfrac{17}{2}\big)
-2\frac{d^2}{dx^2}\log\text{W}\big[\tilde{\phi}_{2}^\ai\big(x;\tfrac{17}{2}\big), \tilde{\phi}_{1}^\ait\big(x;\tfrac{17}{2}\big)\big](x)+4 \nonumber \\
 \qquad {} =x^2+\frac{255}{4x^2}-14-\frac{48}{(6+x^2)^2}+\frac{4}{6+x^2}\\
 \qquad \quad  {}-\frac{640(630+175x^2+12x^4)}{(840+280x^2+30x^4+x^6)^2}+\frac{4(-160+3x^2)}{(840+280x^2+30x^4+x^6)}.
\end{gather*}
As shown clearly above, the singularity at $x^2=-6$
of the new potential corresponds to the characteristic exponents $\rho=-1,2$.

\subsubsection{Global solutions  (F)}\label{sec:sol-F}

The global solutions are singular in $(0,\infty)$
\begin{gather}
\mathcal{H}^F\phi^F_n(x) =\mathcal{E}_n\phi^F_n(x),
\qquad \mathcal{H}^F\eqdef-\frac{d^2}{dx^2}+U_F(x),\qquad \mathcal{E}_n=4n,
\label{schF}\\
\phi^F_n(x) =\frac{e^{-\eta/2}\eta^{-13/4}}{(-6+\eta)^2(-14+\eta)}\mathcal{L}^F_n(\eta),\qquad n=0,1,\ldots,
\label{Fwavfun}
\end{gather}
in which  $\mathcal{L}^F_n(\eta)$ is a degree $n+3$  polynomial in $\eta$ def\/ined by
\begin{gather*}
 \mathcal{L}^F_n(\eta)  \eqdef
 -4\bigl({-}9\big({-}280+140\eta-22\eta^2+\eta^3\big)+n(-14+\eta)(-6+\eta)\bigr)L_n^{(-7)}(\eta)\\
\hphantom{\mathcal{L}^F_n(\eta)  \eqdef}{} -16\eta(-6+\eta)(-12+\eta)
 \partial_\eta L_n^{(-7)}(\eta).
\end{gather*}

They are not square integrable and thus they are not eigenfunctions.
  Some lower members are
\begin{gather}
 \mathcal{L}^F_0(\eta) =36\big({-}280+140\eta-22\eta^2+\eta^3\big),\label{Fgr}\\
\mathcal{L}^F_1(\eta) = -32\big({-}1512+504\eta + 12 \eta^2 -16\eta^3+\eta^4\big),\\
\mathcal{L}^F_2(\eta) =14\big({-}6480+1080\eta+396\eta^2-42\eta^3-12\eta^4+\eta^5\big).\nonumber
\end{gather}

By rewriting the Schr\"odinger equation for $\phi^F_n(x)$ \eqref{schF} in terms of
\begin{equation*}
\phi^F_n(x)=e^{-\eta/2}\eta^{-13/4}y(\eta),
\end{equation*}
we obtain a second-order dif\/ferential equation with a regular singularity at $\eta=0$ and
$\eta=6,14$ and an irregular singularity at $x=\infty$
\begin{gather}
\eta y''-(6+\eta)y'-4\frac{1008-12\eta-16\eta^2+\eta^3}{(-6+\eta)^2(-14+\eta)^2}y+n  y=0,\qquad   n=0,1, \ldots  .
\label{difeqF}
\end{gather}
For each $n$ we have a singular global solution with a $\rho=-2$ and a $\rho=-1$ apparent singularities
\begin{equation*}
y_n(\eta)=\frac{\mathcal{L}^F_n(\eta)}{(-6+\eta)^2(-14+\eta)},\qquad n=0,1,\ldots  .
\end{equation*}
The dif\/ferential equations for the polynomials $\{\mathcal{L}^F_n(\eta)\}$ read
\begin{gather*}
 \eta(-6+\eta)(-14+\eta)\mathcal{L}^F_n(\eta)''-\big(504-104\eta -8\eta^2 + \eta^3\big)\mathcal{L}^F_n(\eta)'\\
\qquad {} +\bigl({-}252-8\eta+3\eta^2+n(-6+\eta)(-14+\eta)\bigr)\mathcal{L}^F_n(\eta)=0,\qquad
n=0,1,\ldots .
\end{gather*}

The  wavefunction $\phi^F_{0}(x)$ \eqref{Fwavfun} with $\mathcal{L}^F_{0}$ \eqref{Fgr}
is related to the potential $U_F(x)$ \eqref{potF} through the prepotential
$w_F(x)$
\begin{gather*}
e^{w_F(x)}\eqdef \frac{e^{-x^2/2}x^{- 13/2 }(-280+140x^2-22x^4+x^6)}{(-6+x^2)^2(-14+x^2)}\\
\hphantom{e^{w_F(x)}\eqdef}{}
\Rightarrow \quad \left(\frac{dw_F(x)}{dx}\right)^2
+\frac{d^2w_F(x)}{dx^2}=U_F(x).
\end{gather*}
By construction, the system is shape-invariant.
The equations and the wavefunctions of~(E) and~(F) systems  have very similar forms.

Here we provide a summary table of (A) to (F), with the seed solution, the potential, the orthogonality weight function and the ``degree'' of the extra member, that is the ``$n$" of the extra member $\mathcal{L}_n$. Case (C) is omitted because it is identical with (B).

\begin{center}
Summary table\\
\begin{tabular}{|c|c|c|c|c|}
\hline
& seed & potential & weight & extra ``deg."\\
\hline
(A) & $e^{x^2}(3+4x^2)^3$ & \tsep{5pt}\bsep{3pt} $x^2-\frac{3}{16x^2}+\frac{48(-3+4x^2)}{(3+4x^2)^2}-\frac{13}{2}$&
$\frac{e^{-\eta}\eta^{1/4}}{(3+4\eta)^4}$&$-2$\\
\hline
(B)& $e^{x^2}(3+4x^2)^3$ & $U_A+1$ & \tsep{5pt}\bsep{3pt}  $\frac{e^{-\eta}\eta^{-1/4}}{(3+4\eta)^4}$&$-3$\\
\hline
(D)&$x^{-3/2}(-3+4x^2)^3$& \tsep{5pt}\bsep{3pt} $x^2-\frac{3}{16x^2}+\frac{48(-9+4x^2)}{(-3+4x^2)^2}-\frac{11}{2}$&
n.a. &$-2$\\
\hline
(E)&$(x^2+6)^3(x^2+14)$& \tsep{5pt}\bsep{3pt} $x^2+\frac{195}{4x^2}+\frac{12(x^2-6)}{(x^2+6)^2}+\frac{4(x^2-14)}{(x^2+14)^2}-16$& $\frac{e^{-\eta}\eta^7}{(6+\eta)^4(14+\eta)^2}$& none\\
\hline
(F)& $(x^2-6)^3(x^2-14)$& \tsep{5pt}\bsep{3pt} $x^2+\frac{195}{4x^2}+\frac{12(x^2+6)}{(x^2-6)^2}+\frac{4(x^2+14)}{(x^2-14)^2}+12$& n.a.& none\\
\hline
\end{tabular}
\end{center}

\subsection{Three-index case}
\label{sec:three-ind}

Among various possibilities of the three index cases with lower degrees, the following has been noted to have a cubic zero:
\begin{equation*}
\text{W}\big[\tilde{\phi}_{1}^\ai\big(x;\tfrac{53}{2}\big),
\tilde{\phi}_{1}^\ait\big(x;\tfrac{53}{2}\big),\tilde{\phi}_{2}^\ait\big(x;\tfrac{53}{2}\big)\big](x)=
-4\frac{e^{-x^2/2}}{x^{51/2}}\big(30+x^2\big)^3\big(390+39x^2+x^4\big).
\end{equation*}
This case is  denoted by  (H).
The potential is
\begin{gather*}
U _H (x) = x^2+\frac{2499}{4 x^2}-52
-\frac{720}{\left(30\!+ \!x^2\right)^2}+\frac{12}{30\!+\! x^2}
 -\frac{312 x^2}{\left( 390\!+\!39x^2\!+\!x^4 \right)^2}+\frac{4 \left(2 x^2-39\right)}{390\!+\!39x^2\!+\!x^4} .
\end{gather*}
The global wavefunction of the Schr\"odinger equation with the eigenvalue $\mathcal{E}_n=4n $ is
\begin{gather*}
\phi^H_n(x)=\frac{e^{-\eta/2}\eta^{51/4}}{(30+\eta )^2(390+39\eta +\eta ^2)}\mathcal{L}^H_n(\eta),\qquad n=0,1,\ldots,
\end{gather*}
in which $\mathcal{L}^H_n(\eta)$ is a degree $n+4$  polynomial in $\eta$.
This belongs to the generic case, and the groundstate corresponds to $\mathcal{L}^H_0 (\eta )$,
which is proportional to $425880 +67704 \eta +4004 \eta^2 +104 \eta^3 +\eta^4 $.
The  polynomial $\mathcal{L}^H_n(\eta) $ also satisf\/ies a second-order dif\/ferential equation with regular singularities at $\eta=0, -30 ,-39/2  \pm \sqrt{-39}/2$ and an irregular singularity at $x=\infty$.
The orthogonality measure for  the polynomials $\{ \mathcal{L}^H_n(\eta) \} $ is given by $\frac{e^{-\eta}\eta^{25}d\eta}{2(30+\eta )^4(390+39\eta +\eta ^2)^2}$.

At present we do not have a systematic means to locate
seed solutions with  high multiplicity of zeros,
we stop our investigation here.

\section{Summary and comments}
\label{summary}

Several families of inf\/initely many global solutions are presented
for certain second-order dif\/fe\-ren\-tial equations \eqref{difeqA},
\eqref{difeqB}, \eqref{difeqD}, \eqref{difeqE} and \eqref{difeqF} having two or more
f\/inite regular singularities and one irregular singularity at inf\/inity.
The characteristic exponent of the global solutions at each f\/inite singularity is $-2$, $-1$ or $0$.
They are obtained from the non-generic as well as
the generic two-indexed Laguerre
polynomials by choosing the parameter $g$ in such a way that
the seed solutions have a triple zero.
In two families (D) and (F),  the extra f\/inite singularity is located inside
the domain of the radial oscillator potential $(0,\infty)$.
In the other three families, the extra singularity is outside of
$(0,\infty)$ and thus they form orthogonal polynomials over
$(0,\infty)$ with the weight functions
(A) $W_{A}(\eta )=e^{-\eta }\eta ^{1/4}/(3+4\eta )^4$,
(B) $W_{B}(\eta )=e^{-\eta }\eta ^{-1/4}/(3+4\eta )^4$,
(E) $W_E(\eta )=e^{-\eta }\eta ^{7}/\{ (6+\eta )^4(14+\eta )^2 \}$.
Since they are obtained within the framework of multi-indexed orthogonal polynomials, their degrees have ``gaps''.
Those polynomials belonging to (A) and (B) were derived in a dif\/ferent way in~\cite{gomez4}.

By construction, these  families of two-indexed Laguerre
polynomials are also the main part of eigenfunctions of
exactly solvable quantum mechanics,
belonging to the class called deformed radial oscillator potentials.
In contrast with the generic multi-indexed polynomials~(E) and~(G),
 the  systems obtained from the non-generic multi-indexed polynomials~(A) and~(B) are not shape invariant.
 By the standard Darboux--Crum transformations, these non-generic two-indexed
 polynomial systems are mapped to the generic one-indexed polynomial systems.

 By using the same method, we tried in vain to f\/ind other global
 solutions with higher degrees ($\rho=-3,-4$) of apparent singularities within
 the multi-indexed Laguerre polynomials.
 It would be interesting to explore the Jacobi polynomials
 counterparts of those results established in the present paper \cite{hst}.
 We wonder if the ideas of the present paper could be generalised to
 a~certain kind of dif\/ference equations in discrete quantum mechanics
 \cite{gos,os15,os8,os13,os24,os22, os12,os7}\footnote{The dual $q$-Meixner polynomial in \cite[\S~5.2.4]{os12} and dual $q$-Charlier
polynomial in \cite[\S~5.2.8]{os12} should be deleted because the hermiticity of
the Hamiltonian is lost for these two cases.}.
 The multi-indexed orthogonal polynomials are also established in
 the framework of dif\/ference Schr\"odinger equations,
 the multi-indexed ($q$-)Racah, Wilson and Askey--Wilson polynomials
 \cite{os20,os23,os17,os26,os27}.

\subsection*{Acknowledgements}

R.S.~is supported in part by Grant-in-Aid for Scientif\/ic Research
from the Ministry of Education, Culture, Sports, Science and Technology
(MEXT), No.23540303 and No.22540186.
K.T.~is supported in part by the Grant-in-Aid for Young Scientists from the Japan Society for the Promotion of Science (JSPS), No.22740107.

\pdfbookmark[1]{References}{ref}
\LastPageEnding

\end{document}